\documentclass[12pt]{article}
\title{Julia directions for
holomorphic curves}
\author{A. Eremenko\thanks{Supported by NSF grant DMS-950036}}
\begin{document}
\date{1988-94}
\maketitle

\begin{abstract} A theorem of Picard type is proved
for entire holomorphic mappings into projective varieties.
This theorem has local nature in the sense that
the existence of Julia directions can be proved
under natural additional assumptions.
An example is given which shows that Borel's theorem
on holomorphic curves omitting hyperplanes has no such
local counterpart.
\end{abstract}

Let ${\bf P}^m$ be complex projective space of dimension $m$ and
$M\subset {\bf P}^m$ be a projective variety. By a divisor on $M$
we mean an intersection of a hyperplane in ${\bf P}^m$ with $M$.
We study holomorphic curves $f:{\bf C}\rightarrow M$.
\vspace{.1in}

\noindent
{\bf Theorem 1.} {\em Every holomorphic map $f:{\bf C}\rightarrow M,$
omitting $2n+1$ divisors such that any $n+1$ of them have empty
intersection, is constant.}
\vspace{.1in}

\noindent
{\em Remark.} The dimension of $M$ is not mentioned in this formulation.
Only the intersection pattern is relevant. 
\vspace{.1in}

\noindent
{\bf Corollary.} {\em Every holomorphic map ${\bf C}\rightarrow {\bf P}^n$,
omitting $2n+1$ hypersurfaces, such that any $n+1$ of them
have empty intersection,
is constant.}  
\vspace{.1in}

\noindent
This Corollary also follows from the results of M. Green \cite{green}
and V. F. Babets \cite{babez}. Their proofs were based on Borel's theorem
(which we will state later). We start with a simple proof of Theorem 1,
independent of Borel's theorem. The method of the proof first appeared
in \cite{ES}. It also provides a new proof of the classical
Picard theorem \cite{lewis,ransford} as well as its generalizations
to quasiregular maps in ${\bf R}^n$ \cite{EL,lewis,HR}.
\vspace{.1in}

{\sl Proof of Theorem 1.} Let $P_1,\ldots,P_{2n+1}$ be the linear forms in
$m+1$ variables defining the divisors. 
Consider a homogeneous representation $F=(f_0:\ldots:f_m)$
of the curve $f,$
where $f_j$ are entire functions without common zeros. Define the subharmonic
function
$$u=\log ||F||=\frac{1}{2}\log (|f_0|^2+\ldots+|f_m|^2).$$ 
Suppose that $f$ is not constant.
Then we may assume that the Riesz
measure\footnote{We call it Cartan measure of $f$. Notice the formula
$T(r,f)=\int_0^r \mu(D(0,t))dt/t$.}
 $\mu$ of $u$ is infinite (if this is not the
case, we can replace $f$ by $f\circ g$ with any transcendental entire $g$).

The functions
$$u_j=\log |P_j\circ F|=\log|P_j(f_0,\ldots,f_m)|,\quad j=1,\ldots 2n+1,$$
are harmonic in ${\bf C}$. 

Let $I\subset \{1,\ldots,2n+1\},\;{\rm card}\, I=n+1$.
Let $\pi:{\bf C}^{m+1}\rightarrow {\bf P}^m$ be the standard
projection.
If $z\in {\bf C}^{m+1},\; ||z||=1$ and $\pi(z)\in M$ then
for some constants $C_1$ and $C_2$ we have
$$C_1\leq \max_{j\in I} |P_j(z)|\leq C_2$$
This follows from the assumption that the intersection of any $n+1$ divisors
is empty. Using the homogeneity we conclude that
$$C_2 ||F(z)||\leq \max_{j\in I} |P_j\circ F(z)|
\leq C_2||F(z)||, \quad z\in {\bf C},$$
so
\begin{equation}
\label{condition}
\max_{j\in I}u_j=u+O(1),\quad {\rm card}\, I=n+1.
\end{equation}
In particular
\begin{equation}
\label{max}
\max_{1\leq j\leq 2n+1} u_j=u+O(1).
\end{equation}
We use the notation $D(a,r)=\{ z\in {\bf C}:|z-a|<r\}$.
\vspace{.1in}

\noindent
{\bf Lemma.} {\em Let $\mu$ be a Borel
measure in ${\bf C},\;\mu({\bf C})=\infty.$
Then there exist sequences $a_k\in{\bf C},\,a_k\rightarrow\infty$ and $r_k>0$
such that
\begin{equation}
\label{infty}
M_k=\mu(D(a_k,r_k))\rightarrow\infty
\end{equation}
and}
\begin{equation}
\label{bounded}
\mu(D(a_k,2r_k))\leq 200\mu(D(a_k,r_k)).
\end{equation}
\vspace{.1in}

\noindent
This Lemma is due to S. Rickman \cite{Rickman}. His formulation contains a
minor mistake (see the discussion below). The Lemma was also used in \cite{EL}.
In the end of the paper we will prove the lemma for completeness.

Apply the Lemma to the Riesz measure $\mu$ of the function $u$.
We obtain two sequences $a_k$ and $r_k$, such that (\ref{infty}) and 
(\ref{bounded}) are satisfied.
Consider the functions defined in $D(0,2)$:
$$u_k(z)=\frac{1}{M_k}(u(a_k+r_k z)-{\tilde u}(a_k+r_kz))$$
and
$$u_{j,k}(z)=\frac{1}{M_k}(u_j(a_k+r_k z)-{\tilde u}(a_k+r_kz)),\quad 1\leq j\leq 2n+1,$$
where ${\tilde u}$ is the smallest harmonic majorant of $u$ in the disc
$D(a_k, 2r_k)$.
The functions $u_k$ are Green potentials that is
$$u_k(z)=-\int_{D(0,2)} G(z,.)\, d\mu_k,$$
where $G(z,.)$ is the Green function of $D(0,2)$ with pole at the point $z$
and $\mu_k$ is the Riesz measure of $u_k$.

It follows from 
(\ref{bounded}) that $\mu_k(D(0,2))\leq 200$ so 
after selecting a subsequence we may assume that $u_k\rightarrow v$,
where $v$ is a subharmonic function, not identically equal to $-\infty$.
(Convergence holds in $L^1_{\rm loc}(D(0,2), dxdy)$, and the Riesz measures
converge weakly, 
see \cite[Theorem 4.1.9]{H}). 
In particular $v$ is {\em not harmonic}
because the Riesz measure of ${\overline D(0,1)}$ is at least one in view
of (\ref{infty}).

All functions $u_{j,k}$ are harmonic 
and bounded from above in view of
(\ref{max}), so we may assume that $u_{j,k}\rightarrow v_j$, each $v_j$
being harmonic or identically equal to $-\infty$ in $D(0,2)$.
From (\ref{condition})
and (\ref{infty}) follows
\begin{equation}
\label{limitcondition}
\max_{j\in I} v_j=v,\quad {\rm card}\, I\geq n+1.
\end{equation}
Thus $v$ is continuous.
For every $I\subset \{1,\ldots,2n+1\}$ of cardinality $n+1$ we consider the set
$E_I=\{ z\in D(0,2):v(z)=v_j(z),\,j\in I\}$. From (\ref{limitcondition})
follows that the union of these sets coinsides with $D(0,2)$.

We conclude that at least one set $E_{I_0}$ has positive area.
By the uniqueness theorem for
harmonic functions all functions $v_j$ for $j\in I_0$ are equal.
Applying (\ref{limitcondition}) to $I_0$ we conclude that $v$ is harmonic.
This is a contradiction which proves the theorem.
\vspace{.2in}

\noindent
{\bf Circles de remplissage and Julia directions}. In Rickman's formulation
of the Lemma there is an extra property
\begin{equation}
\label{rickman}
r_k/|a_k|\rightarrow 0,\quad k\rightarrow \infty.
\end{equation}
The following example shows that this property is not granted in general.
Take 
$$\mu(E)=\int_E\frac{dx dy}{x^2+y^2},\quad E\subset {\bf C}.$$
Then all annuli $\{z:2^m\leq |z|\leq 2^{m+1}\}$ have equal measure and we cannot
find a sequence of discs $D(a_k,r_k)$ satisfying (\ref{infty}) and
(\ref{rickman}). However the following is true.
\vspace{.1in}

\noindent
(*){\em Let $\mu$ be a measure in ${\bf C}$ such that the
sequence
$$A_m=\mu(\{z:t^m\leq |z|\leq t^{m+1}\})$$
is unbounded for some $t>1$. Then there exist sequences $a_k\rightarrow
\infty$ and $r_k>0$ such that the conditions $(\ref{bounded})$,
$(\ref{infty})$ and $(\ref{rickman})$ are satisfied.}
\vspace{.1in}

It is clear that the assumption that $A_m$ is unbounded does not depend
on $t>1.$ To prove (*) we pick a sequence $(m)$ of natural numbers
such that $A_{m}\rightarrow\infty$. There is a covering of the annulus
$\{z:t^m\leq|z|\leq t^{m+1}\}$ by $N_m\rightarrow\infty$ discs such that 
at least one of these discs,
say $D_m$, still has large measure and satisfies
(\ref{rickman}). Now take this disc $D_m$ in place of $D(0,R/4)$ in the
proof of the Lemma (see Appendix at the end of the paper).
This proof will give us a disc $D(a_m,r_m)$
which satisfies (\ref{bounded}), (\ref{infty}) and (\ref{rickman}).

A number $\theta\in [0,2\pi)$ is called a {\em Julia direction} for a
holomorphic
curve $f:{\bf C}\rightarrow M$ if for every
system of divisors $D_1,\ldots,D_q$  such that any $n+1$
of them have empty intersection, and for any $\epsilon>0$ all but
at most $2n$ of these divisors have infinitely many preimages in the
angle $\{z: |\arg z-\theta|<\epsilon\}$.

Using the statement (*) instead the Lemma we prove the following:
\vspace{.1in}

\noindent
{\bf Theorem 2.}
{\em If the Riesz measure $\mu$ corresponding to a holomorphic curve $f$
has the property
$$\limsup_{m\rightarrow\infty}\mu(\{z:2^m\leq|z|\leq 2^{m+1}\}\rightarrow
\infty$$
then $f$ has at least one Julia direction.}

Actually under the assumptions of this theorem there exists a sequence
of discs $D(a_k,r_k)$ satisfying (\ref{rickman}) such that in the union of
these discs all but at most $2n$ divisors have infinitely many
preimages. Any accumulation point of the set $\{\arg a_k\}$ is a
Julia direction. Such sequence of discs is called ``circles de
remplissage''. 

The condition on the measure $\mu$ in Theorem 2 is best possible.
Actually there is an explicit description of all meromorphic functions
(that is holomorphic curves in ${\bf P}^1$) having no Julia directions,
which is due to A. Ostrovski (see, for example
 \cite{Montel}). 

The classical theorem of E. Borel can be formulated in the following way:
\vspace{.1in}

\noindent
{\em Let $f:{\bf C}\rightarrow {\bf P}^n$ be a meromorphic curve with
linearly independent components i. e., whose image is not contained in
a hyperplane. Let $L_1,\ldots,L_{n+2}$ be hyperplanes in general position.
Then $f$ cannot omit} $\cup_j L_j$.
\vspace{.1in}

In the following example a curve $f:{\bf C}\rightarrow {\bf P}^n$ 
with linearly independent components
omits locally $2n$ hyperplanes in general position. That is there exists
a covering of the plane by a finite set of angular sectors such that $f$
omits $2n$ hyperplanes in each sector of the set. So there is no
analogue of Julia directions for Borel's theorem and the estimate $2n$
for the number of exceptional divisors is best possible even in the case
when $M={\bf P}^n$ and $f$ is linearly non-degenerate.
\vspace{.1in}

\noindent
{\bf Example.} For simplicity we construct the example only for $n=2$.
The coordinate representation of $f$ is
$(f_0:f_1:f_2)$ where $f_j(z)=\sin (\varepsilon^j z),\;\varepsilon=
\exp(2\pi i/3),\;j=0,1,2.$ The hyperplanes are defined by the vectors
$$(1,0,1)\,,(1,0,2)\,,(1,1,0)\,,(1,2,0)\,,(0,1,1)\,,(0,1,2).$$
A direct computation (or drawing a picture)
shows that the system is admissible.
In the angle $0<\arg z<\pi/3$, we have $f_2=o(f_0)$ and $f_2=o(f_1)$
as $z\rightarrow\infty$. So the hyperplanes defined by
$(1,0,1)\,,(1,0,2)\,,(0,1,1)$ and $(0,1,2)$, are omitted in this angle.
In the angle
$\pi/4<\arg z<3\pi/4$ we have $f_2=o(f_0)$ and $f_1=o(f_0)$, so the
hyperplanes defined by 
$(1,0,1)\,,(1,0,2)\,,(1,1,0)$ and $(1,2,0)$ are omitted in this angle.
The other angles are 
studied similarly using the property that $f(\varepsilon z)$ is obtained from
$f(z)$ by permutation of coordinates.
\vspace{.1in}

\noindent
{\bf Appendix}.
{\sl Proof of the Lemma.} Take a large number $R>0$, so that $\mu(D(0,R/4))$
is large. Denote $\delta(z)=(R-|z|)/4$. Then find such $a\in D(0,R)$ that
\begin{equation}
\label{lewis}
\mu(D(a,\delta(a))>\frac{1}{2}\sup_{z\in D(0,R)}\mu(D(z,\delta(z))).
\end{equation}
Take $r=\delta(a)$. Then the disc $D(a,2r)$ can be covered by at most
100 discs of the form $D(z,\delta(z))$, so by (\ref{lewis})
$$\mu(D(a,2r))\leq 200\mu(D(a,r)).$$
Putting $z=0$ in (\ref{lewis}) we get
$$\mu(D(a,r))\geq \frac{1}{2}\mu(D(0,R/4)).$$
So we have constructed the disc of arbitrarily large measure and property
(\ref{bounded}). This proves the lemma.

The author thanks Min Ru and Yum-Tong Siu for helpful comments.

{\em Purdue University, West Lafayette, IN 47907 USA}

\begin{thebibliography}{1}
\bibitem[1]{babez} V. F. Babets, Theorems of Picard tipe for 
holomorphic mappings, Siberian Math. J., 25 (1984).
\bibitem[2]{EL} A. Eremenko, J. Lewis, Uniform limits of certain $A$-harmonic
functions with applications to quasiregular mappings, Ann. Acad. Sci. Fenn.,
Ser, A. I. Math., Vol. 16, 1991, 361-375.
\bibitem[3]{ES} A. Eremenko and M. Sodin, The value distribution of
meromorphic functions and meromorphic curves from the point of view of
potential theory, St. Petersburg Math. J., 3 (1992), 109-136.
\bibitem[4]{green} M. Green, Some Picard theorems for holomorphic maps,
Amer. J. Math., 97 (1975), 43-75.
\bibitem[5]{HR} I. Holopainen, S. Rickman, A Picard type theorem for
quasiregular mappings of ${\bf R}^n$ into $n$-manifolds with many ends,
Revista Mat. Iberoamericana, 8 (1992), 131-148.
\bibitem[6]{H} L. H\'{o}rmander, The Analysis of Linear partial Differential
Operators I, Springer, NY, 1983.
\bibitem[7]{lewis} J. Lewis, Picard's theorem and Rickman's theorem by way
of Harnack inequality, Proc. Amer. math. Soc., 122 (1994), 199-206.
\bibitem[8]{Montel} P. Montel, Le\c{c}ons sur les familles normales
de fonctions analytiques et leurs applications,
Paris, Gauthier-Villars, 1927.
\bibitem[9]{ransford} Th. Ransford, Potential theory in the complex plane,
Cambridge Univ. press, Cambridge, 1995.
\bibitem[10]{Rickman} S. Rickman, On the number of omitted values of entire
quasiregular mappings, J. d'Analyse Math., 37, 1980, 100-117.
\end{thebibliography}
\end{document}